%% file: arxiv.tex
\documentclass[12pt]{article}
\RequirePackage{amsthm,amsmath,amsfonts,amssymb,bbm}
\RequirePackage[authoryear]{natbib}
\RequirePackage[colorlinks,citecolor=blue,urlcolor=blue]{hyperref}
\RequirePackage{graphicx}
\usepackage{subfigure}
\usepackage{float}
\usepackage{multirow}
\usepackage{diagbox}
\usepackage[figuresright]{rotating}
\usepackage{pdflscape}
\usepackage{booktabs}
\usepackage{bm}
\theoremstyle{plain}

\newtheorem{theorem}{Theorem}[section]

\newtheorem{remark}{Remark}[section]
\newtheorem{assumption}{Assumption}
\theoremstyle{remark}

\newtheorem{proposition}{Proposition}[section]

\usepackage{tikz}

\newcounter{cnstcnt}

\newcounter{bnstcnt}

\newcounter{dbnstcnt}

\newcounter{dcnstcnt}
\newcommand{\xdc}{%
   \refstepcounter{dcnstcnt}%
   \ensuremath{B_{\thedcnstcnt}}
}
\newcommand{\jdc}[1]{\ensuremath{B_{\ref{#1}}}}

\makeatletter
\newcommand*{\centerfloat}{%
  \parindent \z@
  \leftskip \z@ \@plus 1fil \@minus \marginparwidth
  \rightskip \leftskip
  \parfillskip \z@skip}
\makeatother

\usepackage{enumerate}
\newtheorem{manualtheoreminner}{Assumption}
\newenvironment{manualtheorem}[1]{%
  \renewcommand\themanualtheoreminner{#1}%
  \manualtheoreminner
}{\endmanualtheoreminner}

\newcommand{\x}{X}
\newcommand{\xx}{\boldsymbol{x}}
\newcommand{\zz}{\boldsymbol{z}}
\newcommand{\relu}{\sigma}
\newcommand{\sob}{\mathcal{W}}
\newcommand{\rkhs}{\mathcal{H}_{NT}}
\newcommand{\dom}{\Omega}
\newcommand{\kont}{K^{\mathrm{NT}}}

\newcommand{\kmnt}{K^{\mathrm{NT}}_{\theta}}

\newcommand{\kntt}{K^{\mathrm{NT}}_{t}}

\newcommand{\Pe}{P}
\newcommand{\Pn}{\mathbb{P}_{n}}
\newcommand{\ppp}{\mathbb{P}}
\newcommand{\loss}{l}
\newcommand{\Loss}{L_{n}}
\newcommand{\tJ}{\tilde{J}}
\newcommand{\feta}{f}
\newcommand{\fo}{\feta_{0}}
\newcommand{\beo}{\beta_{0}}
\newcommand{\so}{s_{0}}

\newcommand{\poly}{\operatorname{poly}}
\newcommand{\ddt}{\frac{d}{dt}}
\newcommand{\hft}{\hat{\feta}_{t}}
\newcommand{\hbt}{\hat{\beta}_{t}}
\newcommand{\htt}{\hat{\theta}_{t}}

\newcommand{\tft}{\tilde{\feta}_{t}}
\newcommand{\tbt}{\tilde{\beta}_{t}}
\newcommand{\hfs}{\hat{\feta}_{s}}
\newcommand{\hbs}{\hat{\beta}_{s}}

\newcommand{\tfs}{\tilde{\feta}_{s}}
\newcommand{\tbs}{\tilde{\beta}_{s}}
\newcommand{\tfw}{\tilde{\feta}_{\infty}}
\newcommand{\tbw}{\tilde{\beta}_{\infty}}

\newcommand{\tLoss}{\tilde{L}_{n}}

\newcommand{\Lo}{L_{0}}

\newcommand{\lossbeofo}{\loss_{\beo, \fo}}
\newcommand{\lossbf}{\loss_{\beta, \feta}}
\newcommand{\ddh}{\mathbf{h}}
\newcommand{\seff}{\tilde{S}}
\newcommand{\tddh}{\tilde{\ddh}}
\newcommand{\thh}{\tilde{h}}
\newcommand{\idicator}{\mathbbm{1}}

\newcommand{\stopt}{t_{s}}
\newcommand{\hbts}{\hat{\beta}_{\stopt}}
\newcommand{\hfts}{\hat{\feta}_{\stopt}}

\usepackage{newtxtext} 
\usepackage{libertine}  
\begin{document}

\def\spacingset#1{\renewcommand{\baselinestretch}%
{#1}\small\normalsize} \spacingset{1}


\title{\bf{Semiparametric M-estimation with overparameterized neural networks}}
	
 \author{Shunxing Yan \thanks{E-mail: sxyan@stu.pku.edu.cn},~~~Ziyuan Chen\thanks{E-mail: chenziyuan@pku.edu.cn}~~\mbox{ and }~Fang Yao\thanks{		Fang Yao is the corresponding author, E-mail: fyao@math.pku.edu.cn. } \\
 \hspace{.2cm} 
 School of Mathematical Sciences, Center for Statistical Science, \\	Peking University, Beijing, China\\
 }
\date{}
\maketitle
\vspace{-0.32in} 
\input{chapters/abstract}

\spacingset{1.39} 

\input{chapters/mainpart}

\bibliographystyle{apalike}
\bibliography{bibliography}  

\end{document}

%% file: chapters/abstract.tex
\begin{abstract} 

We focus on semiparametric regression that has played a central role in statistics, and exploit the powerful learning ability of deep neural networks (DNNs) while enabling statistical inference on parameters of interest that offers interpretability. Despite the success of classical semiparametric method/theory, establishing the $\sqrt{n}$-consistency and asymptotic normality of the finite-dimensional parameter estimator in this context remains challenging, mainly due to nonlinearity and potential tangent space degeneration in DNNs. In this work, we introduce a foundational framework for semiparametric $M$-estimation, leveraging the approximation ability of overparameterized neural networks that circumvent tangent degeneration and align better with training practice nowadays. The optimization properties of general loss functions are analyzed, and the global convergence is guaranteed. Instead of studying the ``ideal'' solution to minimization of an objective function in most literature, we analyze the statistical properties of algorithmic estimators, and establish nonparametric convergence and parametric asymptotic normality for a broad class of loss functions. These results hold without assuming the boundedness of the network output and even when the true function lies outside the specified function space. 
To illustrate the applicability of the framework, we also provide examples from regression and classification, and the numerical experiments provide empirical support to the theoretical findings.

\end{abstract}

%% file: chapters/mainpart.tex
\section{Introduction}\label{sec:intro}
\subsection{Literature review}
Semiparametric models constitute an essential category of statistical methods, involving both finite and infinite-dimensional parameters. 
Typically, more consideration is given to the finite-dimensional components and the latter is regarded as a nuisance \citep{van2000asymptotic, kosorok2008introduction}. 
Compared to parametric models with only finite-dimensional parameters, semiparametric models provide stronger representation capabilities. Moreover, in contrast to nonparametric models, they alleviate the curse of dimensionality and allow easier inference on the finite-dimensional parameters of primary interest. 
Benefiting from these advantages, semiparametric models have garnered increasing attention in statistics over an extended period, including regression analysis \citep{
ahmad2005efficient,liang2010estimation}, 
survival analysis \citep{cox1972regression, cox1975partial, huang1999efficient,zeng2007maximum}, causal inference \citep{rosenbaum1983central, chernozhukov2018double} 
and many others; 
see \citet{
tsiatis2006semiparametric,kosorok2008introduction,horowitz2012semiparametric} 
for a comprehensive introduction. 
Among these works, one of the most common estimation methods is the semiparametric $M$-estimation \citep{geer2000empirical,ma2005robust,cheng2010bootstrap}, where the estimators for both parametric and nonparametric parameters are obtained simultaneously by minimizing (or maximizing) certain objective functions.
There has been much research on this method, employing classical nonparametric regression techniques, including linear sieves \citep{ding2011sieve,ma2021estimation,tang2023survival},
local polynomial method with closed-form representations \citep{wang2010estimation,liang2010estimation} 
and penalized methods \citep{mammen1997penalized,ma2005penalized}.
In general, efficient asymptotic $\sqrt{n}$-normality of finite-dimensional parameters and minimax optimal convergence of the nonparametric parts can be established. 
However, most existing works based on traditional techniques often face challenges when applied to complex structured data increasingly encountered in modern applications.
Consequently, there is growing interest in developing estimators based on advanced methodologies such as neural networks, with rigorous theoretical guarantees, motivated by the significant achievements of deep learning in recent years.

Deep neural networks (DNNs), especially those with ReLU activation function, have received significant attention in recent studies due to their strong learning abilities. 
Approximation bounds for ReLU DNNs have been established in various function settings \citep{yarotsky2017error,yarotsky2018optimal, schmidtaos,suzuki2018adaptivity,yarotsky2020phase,kohler2021rate}, which are shown to be (nearly) optimal in terms of both width and depth \citep{shen2020,lu2021deep,shen2022optimal}. 
A key advantage of DNN estimators, distinguishing them from classical nonparametric regression estimators, is their adaptivity to various low intrinsic-dimensional structures, such as the hierarchical composition model where the target function follows a specific compositional structure \citep{bauer2019deep, schmidtaos, kohler2021rate}, and cases with low-dimensional inputs  \citep{schmidt2019deep,chen2019nonparametric,cloninger2021deep,nakada2020adaptive,huangjian}. 
Therefore, DNNs are increasingly recognized as an important nonparametric approach in a wide range of statistical problems, including nonparametric regression \citep{bauer2019deep,schmidtaos,suzuki2018adaptivity,kohler2021rate, wang2024deep}, survival analysis \citep{deephazard,deepcox}, causal inference \citep{farrell2021deep,chen2024causal}, factor augmented and interaction models \citep{fan2022factor,bhattacharya2023deep},
repeated measurements model \citep{yan2023nonparametricregressionrepeatedmeasurements}, among others.

As a common ground, most of the previously mentioned statistical works study the generalization ability of the (nearly) empirical risk minimizers and bound the estimation error by network size relative to the training sample. 
However, the optimization landscape in deep learning presents unprecedented difficulty due to nonlinearity and nonconvexity, leading to the estimation more likely being local minimizers. 
Moreover, in practice, neural networks are often trained with overparameterization, i.e., the network parameters may vastly exceed the training sample size, while still avoiding the traditional pitfall of overfitting \citep{zhang2021understanding}. This does not align with the statistical analysis in the previously mentioned works. 
Fortunately, there have been meaningful results to address these gaps. 
For instance, \citet{jacot2018neural} introduced a framework that analyzes the training of wide neural networks, drawing an analogy to kernel regression. 
Specifically, they compared the gradient flow of the least squares loss with that of kernel regression and demonstrated that, as the network width approaches infinity, the Neural Tangent Kernel (NTK) converges to a time-invariant limit \citep{arora2019exact,lee2019wide,bietti2019inductive,geifman2020similarity,chen2020deep,bietti2020deep,lai2023generalization}. 
From the perspective of optimization, wide neural networks with random initialization have been shown, with high probability, can achieve global minimization through gradient-based methods \citep{DBLP:journals/corr/abs-1810-02054, li2018learning, allen2019convergence}. 
In terms of statistical generalization performance, \citet{hu2021regularization,suh2021non} established the convergence rate of penalized least square regression problem, and more recent studies have examined least squares regression with early stopping, providing uniform convergence guarantees for neural network kernels \citep{lai2023generalization,li2023statistical,lai2023residual}.

\subsection{Challenges and main contributions}\label{sec:chacon}
Despite the impressive empirical performance of deep learning models, they are commonly regarded as black boxes, lacking interpretability and theoretical support. Semiparametric modeling provides a useful way for making interpretable inferences while leveraging the learning ability of neural networks. 
Training DNNs involves a large-scale network parameter learning via loss minimization, and aligns naturally with the framework of semiparametric $M$-estimation which simultaneously estimates nonparametric (via network weights) and parametric parameters. 
Several studies have explored this problem \citep[e.g.,][]{deepcox},  but defaulted to some ``good'' assumptions on the tangent space—a critical issue we shall discuss below. 
Another commonly used method is $Z$-estimation, where the finite-dimensional parameter is typically taken as a functional of the nonparametric component. 
It usually takes two steps: first estimating the nonparametric component and then substituting it into an equation to solve for the parametric component. 
For such two-step procedures, the doubly debiased/robust methods \citep{chernozhukov2018double,farrell2021deep} are employed to achieve efficiency. By comparison, semiparametric $M$-estimation is computationally straightforward via network training, and requires little statistical manifest with broader applicability. While a simpler two-step plug-in method can also attain efficiency \citep{chen2024causal}, one has to impose assumptions on the tangent approximation ability in a similar manner to \citet{deepcox}.

The above discussion inspires our investigation: what foundation of semiparametric $M$-estimation theory is undermined in neural network based models, and how to rebuild it? 
Denote the semiparametric model as $P_{\beta, \feta }$, where \( \beta \in \mathbb{R}^p \) is the finite-dimensional parameter of interest, and \( \feta : \mathbb{R}^d \to \mathbb{R} \) is a nuisance function belonging to a infinite-dimensional space.
The semiparametric $M$-estimation procedure aims to optimize an objective function, expressed as $\Pn \loss_{\beta, \feta}$ with empirical measure integration $\Pn$, either in a suitable sieve or with an additional penalty term. Notably, the information of the parameters, \( \beta \) and \( \feta \), is coupled. 
Since the nuisance parameter \(\feta\) resides in an infinite-dimensional space, it significantly increases the complexity of the estimation problem. Consequently, it is reasonable that the full set of parameters achieves a nonparametric convergence rate that is slower than $n^{-1}$.

A central challenge of semiparametric $M$-estimations is to establish \( \sqrt{n} \)-consistency and the asymptotic normality of the estimator for the finite-dimensional component \( \beta \). 
To decouple the two parts of parameters, taking likelihood estimation as an example, the efficient score is commonly utilized, which removes the effect of nuisance parameters.
Specifically, we define $\Lambda_{\beta,\feta}$ as the tangent set for the nuisance parameter, i.e. it contains score functions $\partial_{\feta} \loss_{\beta,\feta}[h]$ with all appropriate directions $h$, see Section \ref{sec:str} for a rigorous definition. 
Additionally, denote $\Pi_{\theta, \feta}$ as the orthogonal projection operator onto the closure of the linear span of $\Lambda_{\beta,\feta}$. 
Then the efficient score can be constructed by $\partial_{\beta} \loss_{\beta,\feta}$ subtracting its projection on nuisance tangent space $\Lambda_{\beta,\feta}$, i.e. 
$\seff=S_1(\beo, \fo)-\Pi_{\beo, \fo} S_1(\beo, \fo), $
where 
$S_1(\beta, \feta) = 
\partial \lossbf/\partial \beta $.
Accordingly, $\seff$ is orthogonal to the tangent space $\Lambda_{\beta,\feta}$, thus mitigating the nuisance information.
Once the empirical efficient score $\Pn \seff(\hat{\beta}, \hat{\feta})$ is proven to be small enough as $o_p(n^{-1/2})$, asymptotic normality of $\sqrt{n}(\hat{\beta} - \beo)$ can then be established via standard Taylor expansion and entropy analysis \citep{van2000asymptotic,kosorok2008introduction}.

However, establishing that 
$$\Pn \seff(\hat{\beta}, \hat{\feta})  = \Pn  (S_{1}(\hat{\beta}, \hat{\feta}) - \Pi_{\beo, \fo} S_1(\beo, \fo))$$ is $ o_p(n^{-1/2})$ for neural network estimators is highly nontrivial. 
The first term on the right-hand side is typically easy to bound, as $\hat{\beta}$ is a finite-dimensional (nearly) minimizer. The main challenge arises from bounding the second term. 
Because $\Pi_{\beo, \fo} S_1(\beo, \fo)$ is in the the closure of the linear span of $\Lambda_{\beta,\feta}$, there could be some direction $\tddh$ such that $\Pi_{\beo, \fo} S_1(\beo, \fo) = S_{2}\left(\beo, \fo\right)[\tddh ]$. 
Denote the neural network function as $\feta_{\theta}$, where $\theta$ represents the network parameters. 
When $\hat{\feta}_{\theta}$ is an exact local or global minimizer, we can only obtain that $\partial_{\theta} \Pn \loss_{\hat{\beta}, \hat{\feta}_{\theta}} = 0$, but not necessarily $\Pn S_{2}(\hat{\beta}, \hat{\feta})[\tddh ] = 0$. 
Denote the tangent space of the network at $\theta$ by $$\mathcal{T}_{\theta}\mathcal{F}_{\text{NN}} := \operatorname{Span}\{\nabla_{\theta} f_{\theta}(\xx) ^{T}  (\theta_{1} - \theta) : \text{ all } \theta_{1} \text{ having the same shape with }  \theta \}.$$ 
Then, for any $h \in \mathcal{T}_{\hat{\theta}}$, we have $\Pn (S_{2}(\hat{\beta}, \hat{\feta})[ h ]) = 0$. If the components of $\tddh$ lie in, or can be well approximated by, the tangent space $\mathcal{T}_{\hat{\theta}}$, we could conclude that $\Pn S_{2}(\hat{\beta}, \hat{\feta})[\tddh ]$ is sufficiently small. 
Therefore, it is preferable for the tangent space to be large to contain or approximate as many possible directions as possible. 
For traditional linear estimators using a sieve space $\mathcal{S}=\left\{f_\alpha(x): \ f_\alpha(\cdot)=\sum_{i=1}^m \alpha_i b_i(\cdot)\right\}$ with basis functions $\{b_i(\cdot),i=1,\cdots,m\}$ (e.g., regression spline estimators), a large tangent space is straightforward to guarantee, since
the linearity ensures that the tangent space of $\mathcal{S}$ at any $f_\alpha(\cdot)\in \mathcal{S}$ is $\mathcal{T}_\alpha\mathcal{S} = \operatorname{Span} \left\{\partial_{\alpha_i} f_\alpha(\cdot)\right\} \cong \mathcal{S}$. This provides a strong approximation ability. However, this does not necessarily hold for the tangent space of neural networks.

Therefore, the crucial question arises on the approximation ability of the neural network tangent space, and should not be simply suppressed by presupposed conditions as in previous works \citep{deepcox,chen2024causal}:
\begin{quote}
\normalsize
\emph{Whether the tangent space $\mathcal{T}_{\theta} \mathcal{F}_{\text{NN}}$ of a neural network always has good approximation ability at every possible $\theta$? If not, can the semiparametric DNN $M$-estimators still achieve a $\sqrt{n}$-consistency for the finite-dimensional parameters? }
\end{quote}
Unfortunately, the answer to the former question is negative. 
As a toy example, we consider a fully connected neural network in which the weights and biases within each layer take the same values.
Due to the lack of identifiability among network parameters, the derivatives with respect to parameters in equivalent positions are the same. 
Consequently, the dimension of the tangent space $\mathcal{T}_{\theta}\mathcal{F}_{\text{NN}}$ is only linearly related to the depth of the neural network, which is much smaller than the number of parameters (i.e., the square of the width multiplied by the depth).
In such cases, the tangent space $\mathcal{T}_{\theta}\mathcal{F}_{\text{NN}}$ will not provide a sufficiently rich approximation ability for possible $\tddh$. 
As discussed above, whether the efficient score is sufficiently small to establish the $\sqrt{n}$-consistency of the parametric component of the DNN estimator remains an open problem.

Nonetheless, the counterexample presented above is so extreme that it is reasonable to presume that it has only a small probability of occurring in learning. Hence, a workable theoretical treatment requires a more careful analysis of the estimation procedure and randomness introduced by the neural network, like the random initialization. 
Additionally, we hope to address the limitations in most current statistical analyses of deep neural networks: they ignore the nonlinearity and nonconvexity of the optimization landscape and directly consider the ``ideal'' global minimizer. 
Motivated by these considerations, we study overparameterized neural networks that provide a meaningful way for analyzing the optimization and statistical performance of the algorithmic solutions. 
Further, for generality, we consider a broader class of loss functions beyond the least squares criterion, which introduces additional difficulties for both optimization convergence analysis and statistical inference of the algorithmic solution. 
In summary, our main contributions to tackling the challenges are elaborated below.

\begin{itemize}
\item 
Methodologically, for the semiparametric $M$-estimation problem, we employ a new overparameterized neural network estimator, which better aligns with practical settings and facilitates the optimization analysis. The \( l^2 \) penalization of the neural network parameters is applied for regularization, allowing subsequent analyses of general loss rather than only least square regression. This also ensures that the scores at the estimation remain sufficiently small, contributing to the asymptotic normality of the finite-dimensional parameter estimation. 
Then we consider a continuous gradient flow framework to investigate the training dynamics of the neural network. By incorporating random initialization into the analysis, with high probability, we bound the difference between the overparameterized neural network training flow and the ideal RKHS optimization process. 
It also suggests that the counterexamples mentioned above with degenerate tangent space would be only taken with a small probability. 
Furthermore, we establish the global algorithmic convergence of the proposed overparameterized neural semiparametric $M$-estimator, providing a theoretical guarantee of the training procedure.

\item  
Theoretically, we analyze the statistical properties of the algorithmic solution, including the generalization error of the nonparametric component and the $\sqrt{n}$-consistency/asymptotic normality of the parametric component.
Specifically, unlike the common convergence analysis of neural network estimation, which often assumes that the network output is bounded, 
the optimization solution is analyzed directly.
Especially when the true function does not lie in the desired space and the loss function is of a general form instead of least squares, estimators that do not impose bounded assumptions on the nonparametric part present significant challenges to theoretical analysis. To establish the convergence rate, we introduce a new condition, referred to as the ``Huberized margin condition'', which relaxes the standard assumptions and is easier to satisfy by unbounded nonparametric candidate function classes. 
Building on the above results and some regularity conditions, we show that the parametric component achieves $\sqrt{n}$-consistency and asymptotic normality. 
The latter result demonstrates the efficiency for the least favorable submodels and enables interpretable inference for parameters of interest. 
Lastly, we discuss two commonly used models: partially linear regression and classification. In these examples, the aforementioned properties of the proposed estimator are examined to be valid. 
Moreover, the numerical results corroborate our theoretical analysis by the finite sample behavior.

\end{itemize}

\subsection{Notations and organization}

In this paper, we use the notation \( a_n \lesssim b_n \) to indicate that \( a_n \leq C b_n \) for some constant \( C > 0 \) independent of \( n \), with \( \gtrsim \) defined analogously. Then \( a_n \asymp b_n \) means that both \( a_n \lesssim b_n \) and \( b_n \lesssim a_n \). Denote \( a_n = o(b_n) \) if \(a_n/b_n \to 0 \), \( a_n = O(b_n) \) if \( a_n \lesssim b_n \), and \( a_n = \Theta(b_n) \) if \( a_n \asymp b_n \). 
We also use \( a_n = o_p(b_n) \) to indicate that \( a_n/b_n \to_p 0 \), and \( a_n = O_p(b_n) \) if, for any \( \epsilon > 0 \), there exists a constant \( C_\epsilon > 0 \) such that \( \sup_n P\left( \| X_n \| \geq C_\epsilon |Y_n| \right) < \epsilon \).
Additionally, $\poly(a,b,...)$ indicates some polynomial about $(a,b,...)$.
The notation \( \mathbb{I} \) denotes the indicator function. 
For probabilistic integrals, \( \mathbb{P} \) represents the theoretical expectation with respect to the population distribution, while \( \mathbb{P}_n \) denotes the empirical expectation derived from the finite sample.

The rest of the article is organized as follows. 
In Section \ref{sec:pre}, following an overview of overparameterized neural networks and neural tangent kernels, we introduce the semiparametric model framework for neural $M$-estimation trained by the gradient flow algorithm.  Section \ref{sec:str} develops the statistical theory for the algorithmic estimators, addressing the nonparametric convergence of the nonparametric component and the asymptotic normality of the parametric component. In Section \ref{sec:example}, we examine two illustrative examples, the partially linear regression and classification, demonstrating the validity of the theoretical results. 
Finally, the numerical experiments are presented in Section \ref{sec:num}, providing empirical evidence of the proposed method and theory. 

\section{Semiparametric M-estimation with neural networks}\label{sec:pre}
In this section, we first introduce the overparameterized neural network used and some important related concepts. Then we present the considered semiparametric model and neural $M$-estimation framework. Furthermore, a detailed analysis of the optimization convergence would also be provided.

\subsection{Overparameterized neural network and neural tangent kernel}\label{sec:space}
In this paper, we primarily consider feedforward fully connected neural networks. 
Given positive integers $m_0, m_1, ..., m_{L+1}$ with $m_0 = d, m_{L+1}=1$, the network is defined as following 
\begin{equation}\label{equ:fnn0}
\tilde{f}_{\theta}(\xx)=\mathcal{L}_{L} \circ  \relu  \circ \mathcal{L}_{L-1} \circ  \relu  \circ \mathcal{L}_{L-2} \circ  \relu  \circ \cdots \circ \mathcal{L}_{1} \circ  \relu  \circ \mathcal{L}_{0}(\xx), 
\end{equation}
where $\relu:\mathbb{R}^{m_{i}}\rightarrow \mathbb{R}^{m_{i}}$ is the activation function. 
For $i\in\{0,1,...,L-1\}$, $\mathcal{L}_{i}(\boldsymbol{y})=\sqrt{2}(W_{i} \boldsymbol{y}+b_{i}) / \sqrt{m_{i+1}}$,  and for $i=L$, $\mathcal{L}_{i}(\boldsymbol{y})=W_{i} \boldsymbol{y}+b_{i}$, where $W_{i} \in \mathbb{R}^{m_{i+1} \times m_{i}}, b_{i} \in \mathbb{R}^{m_{i}}$ for $\boldsymbol{y}\in \mathbb{R}^{m_i}$. 
In this work, we default $\relu$ to the rectified linear unit (ReLU) function $a\mapsto \max \{a, 0\}$, and most results presented can be generalized to other activation functions. 
Additionally, we use $\theta$ to denote all parameters in the neural network. 
To conveniently analyze the effect of the width $m_i$, we assume that $ m \leq \min_{1\leq i\leq L}  m_i \leq \max_{1\leq i\leq L}  m_i \leq C m$ always holds for some constant $C$.

In practice, the neural networks are overparameterized, which means the width $m$ will be much larger than the sample size $n$. 
To better align with practical settings and facilitate the analysis of optimization properties, we also consider overparameterized networks in the following. 
Before training, random initialization of weights in networks is usually employed to break symmetry and prevent all neurons from learning identical features, thereby improving learning efficiency and promoting better convergence. 
Therefore, we randomly initialize the neural network weight matrices $W_i$'s and the bias $b_0$ from an i.i.d. standard normal distribution, while the other biases $b_i$'s for $i \geq 1$ are initialized to zero. 
Denoting the initial parameters by $\theta_0$, we define
$$
f_{\theta}(\xx) = \tilde{f}_{\theta}(\xx) - \tilde{f}_{\theta_{0}}(\xx)
$$
to ensure that $f_{\theta}$ is initially zero in the training, i.e. $f_{\theta_{0}}(\xx) \equiv 0$, for theoretical convenience.

In the analysis of overparameterized neural networks, the neural tangent kernel (NTK) \citep{jacot2018neural,arora2019exact} is commonly employed.
We now introduce the NTK for finite $m$ as 
$$
   \kmnt (\xx,\xx') = \nabla_{\theta} f_{\theta}(\xx)  ^{T} \nabla_{\theta} f_{\theta}(\xx') .
$$
Fixing $L$ and letting $m \to \infty$, the kernel function converges to 
$$
\kont \left(\xx, \xx'\right) = \lim _{m \rightarrow \infty} \kont_{\theta_{0}}\left(\xx, \xx'\right), 
$$
Recent works have established that this convergence holds uniformly \citep{lai2023generalization}. 
In the sequel, we will generally refer to the limiting kernel as NTK for brevity. 
Moreover, an explicit formula for NTK with the aforementioned random initialization can be derived. 
\begin{proposition}\label{prop:explicit}
Under the random initialization mechanism proposed for $W_i$'s and $b_i$'s, the neural tangent kernel has the explicit expression as 
\begin{equation}\label{equ:101}
    \kont (\xx, \xx' ) = \sum_{l=0}^{L} \Big( \sqrt{ (\|\xx\|_2^2+1) (\|\xx ' \|_2^2+1) }
    \kappa_1^{(l)} \left( \tilde{u} \right) + \idicator_{l\geq 1} \Big) \prod_{r=l}^{L-1}\kappa_0(\kappa_1^{(r)}(\tilde{u})),
\end{equation}
where 
 $
\tilde{u} := 
({\xx}^{T}{\xx}' + 1)/\sqrt{ (\|\xx\|_2^2+1) (\|\xx ' \|_2^2+1) }$, 
$\kappa_0(t)$ and $\kappa_1(t)$ are arc-cosine kernels of degree $0$ and $1$, i.e. 
$$
\kappa_0(t)=\frac{1}{\pi}(\pi-\arccos t), \quad \kappa_1(t)=\frac{1}{\pi}\left[\sqrt{1-t^2}+t(\pi-\arccos t)\right], 
$$
and $h^{(r)}, r \geq 1$ denote the $r$-times composition of a function $h$ 
while $h^{(0)}$ is the identity map. 
\end{proposition}

These expressions presented here differ slightly from previous results \citep{bietti2019inductive,geifman2020similarity}, due to the special initialization of the bias, which simplifies subsequent derivations related to the RKHS of $\kont$ over a general domain. Specifically, we introduce the following equivalence property that facilitates subsequent statistical analysis.

\begin{proposition}\label{prop:p1}
    For any $\dom \subset  \mathbb{R}^d$ with Lipschitz boundary, the RKHS $\rkhs$ associated to $\kont$ is norm-equivalent to the Sobolev space $\sob^{(d+1)/2,2}(\dom)$. 
\end{proposition}

\begin{remark}
Since the properties of Sobolev spaces are extensively studied, given the above proposition, many results can be derived more easily using existing results of Sobolev spaces. 
We emphasize that the smoothness index $(d+1)/2$ is determined by the ReLU activation function. 
Generally, smoother activation functions lead to higher smoothness in the corresponding NTK Sobolev space. 
For example, the rectified power unit activation $a \mapsto \max \{a, 0\}^r$ with positive integer $r$ \citep{DBLP:journals/corr/Bach14, DBLP:journals/corr/abs-2109-06099}, which is weakly differentiable of order $r$, can be proven to lead to Sobolev space $\sob^{(d+2r-1)/2,2}$ via similar analysis as Proposition \ref{prop:p1}. 
\end{remark}

We close this subsection by recalling the motivations behind using overparameterized DNNs in this paper. 
In practical applications, neural networks are often overparameterized, where the number of parameters exceeds the sample size. 
Furthermore, rather than studying the ideal global minimizers as in common statistical works, we hope to investigate the statistical properties of the algorithmic estimators, particularly under the overparameterized optimization convergence guarantee. 
By introducing stochastic initialization and analyzing the optimization process of the overparameterized DNNs, we can avoid, with high probability, the counterexample mentioned in Section \ref{sec:chacon} on the degeneration of the DNN tangent space, and hence helps to establish $\sqrt{n}$-consistency/asymptotic normality of the estimation of the parametric component.

\subsection{Semiparametric neural M-estimation}\label{sec:estimat}

Consider a general semiparametric statistical model $P_{\beta, \feta }$, where $\beta \in \mathbb{R}^{p}$ is a Euclidean parameter of interest and $ \feta : \mathbb{R}^{d} \mapsto \mathbb{R}$ is the nuisance function in an infinite-dimensional space. 
Let $(Y,Z,\x)$ follows the distribution $P_{\beta, \feta }$, where  $Z \in \mathbb{R}^{p}$ and $\x \in \mathbb{R}^{d}$ are finite-dimensional covariates. 
For simplicity, the domains of $Z$ and $\x$ are assumed to be bounded, compact sets with regular boundaries, and their densities are bounded away from zero and infinity. Moreover, we assume that $\mathbb{E}\left[ Y^2|Z,\x \right]$ is finite. 
Consider a nonnegative loss function $\loss_{\beta, \feta} = \loss(Z^{T}\beta, \feta(\x),Y)$ such that the true parameters $(\beo, \fo)$ minimizes the risk $\Pe \loss_{\beta, \feta} = \int \loss(Z^{T}\beta, \feta (\x), Y) d P_{\beta, \feta }$.  
Common choices for $\loss$ include the negative log-likelihood function, squared loss function, or other robust loss functions.

Given i.i.d. observations \( \{(Y_i, Z_i, \x_i), i = 1, 2, \dots, n\} \), we aim to estimate the unknown parameters \( (\beta, \feta) \) by minimizing the criterion
\[
\Loss(\beta, \feta) := \Pn \loss_{\beta, \feta} + \lambda_n J(\feta),
\]
where \( \lambda_n \geq 0 \) is a tuning parameter and \( J \) is a penalty term, within a suitably chosen parameter set \( (\mathbb{R}^p, \mathcal{F}_n) \). Letting the nonparametric function set \(\mathcal{F}_n\) be the overparmeterized neural network class, for each \( \feta_\theta \in \mathcal{F}_n \) with network parameter \( \theta \) and initial parameter \( \theta_0 \), we define the penalty term as
\[
J(\feta_\theta) = \|\theta - \theta_0\|_2^2,
\]
which regularizes the complexity of the neural network to prevent overfitting. 
An alternative regularization technique is early stopping, while existing statistical analyses for early stopping in kernel gradient algorithms mainly focus on least square losses \citep{yao2007early,raskutti2014early}. 
For generality, this work considers a broader class of loss functions, hence, we adopt the penalization strategy.
Even under the penalization framework, there are still considerable challenges in the statistical theory. 
Unlike typical convergence analyses on neural network estimation, which often assume that the nonparametric function is bounded or that the loss function is Lipschitz continuous, we will study the optimization solution without boundedness assumptions. This brings difficulties when the true function does not belong to the desired space and the loss function takes a more general form than the least squares.

\subsection{Gradient flow optimization and convergence}

Given the above considerations, we aim for our estimator to satisfy
\begin{equation}\label{equ:esti}
(\hat \beta, \hat f) = (\hat \beta, f_{\hat \theta}) \quad \text{where} \quad (\hat \beta, \hat \theta) \approx \underset{\beta \in \mathbb{R}^p, f_\theta \in \mathcal{F}_n}{\arg \min} \Loss(\beta, \feta_\theta).
\end{equation}
To optimize the objective in \eqref{equ:esti}, gradient-based algorithms are widely employed, with various modifications such as stochastic sampling and momentum methods being particularly common in practice. A substantial body of literature addresses optimization algorithms for neural networks, highlighting diverse methodologies and their relative effectiveness in improving training outcomes. In this study, for theoretical convenience, we focus on gradient flow, which serves as the continuous counterpart of gradient descent.
Let $\htt$ denote the neural network parameters at time $t \geq 0$. Correspondingly, let $\hft = \feta_{\htt}$ represent the neural network output, and $\kntt=K^{\mathrm{NT}}_{\htt}$ denote the neural tangent kernel at $t \geq 0$. With the initial values set as $\hat{\theta}_0=\theta_0$ and $\hat{\feta}_0(\xx) = 0,$ for all $\xx \in \dom$, the gradient flow training process of the parameters \( (\hbt,\htt)\) is governed by the following equations: 
\begin{equation}\label{equ:011}
    \ddt\hbt=-\nabla_{\beta} 
\Loss \left(\hbt,\htt\right)=-\frac{1}{n}\sum_{i=1}^{n} \loss_{1}'(Z_{i}^{T}\hbt, \hft (\x_{i}), Y_{i}) Z_{i},
\end{equation}
and
\begin{equation}\label{equ:012}
\ddt \htt=-\nabla_{\theta} 
\Loss \left(\hbt,\htt\right)=-\frac{1}{n}\sum_{i=1}^{n} \loss_{2}' (Z_{i}^{T}\hbt, \hft (\x_{i}), Y_{i}) \nabla_{\theta}\hft (\x_{i}) - 2 \lambda_{n} (\htt - \theta_{0}).
\end{equation}
Then the flow of $\feta$ is defined by a dynamical system as 
\begin{equation}\label{equ:013}
\begin{aligned}
\ddt \hft (\xx)& = \nabla_{\theta}\hft (\xx)^{T} \ddt \htt \\
& = -\frac{1}{n}\sum_{i=1}^{n} \loss_{2}' (Z_{i}^{T}\hbt, \hft (\x_{i}), Y_{i}) \kntt (\xx,\x_{i}) - 2 \lambda_{n} \nabla_{\theta}\hft (\xx)^{T} (\htt - \theta_{0}). 
\end{aligned}
\end{equation}

For theoretical analysis, we take an ideal estimator in the RKHS associated with reproducing kernel $\kont$ as a surrogate, which is shown to be sufficiently close to the neural algorithmic estimator \( (\hbt,\htt) \).
In this context, the penalty term is defined as
$
\tJ(\feta ) = \|\feta\|_{\rkhs}^2
$, which is standard and encourages smoother solutions by penalizing the complexity of the function $\feta$ in the RKHS norm. This penalty term leads to the following regularized optimization criterion:
\begin{equation}\label{equ:rkhs}
    \tLoss(\beta, \feta ) = \Pn \loss_{\beta, \feta} + \lambda_{n} \tJ(\feta ),
\end{equation}
where $\Pn \loss_{\beta, \feta}$ still represents the empirical risk term and $\lambda_{n} \tJ(\feta)$ is the tuning parameter.
Consequently, within the RKHS, a gradient flow training process $\{\tbt,\tft\}$ with initial value $\tilde{\theta}_0=0$ and $\tilde{\feta}_0(\xx) = 0, \forall \xx \in \dom$ is adopted as  
$$
\ddt\tbt=-\nabla_{\beta} 
\tLoss \left(\tbt,\tft \right)=-\frac{1}{n}\sum_{i=1}^{n} \loss_{1}'(Z_{i}^{T}\tbt, \tft (\x_{i}), Y_{i}) Z_{i},
$$ 
and 
$$ 
\ddt \tft  = -\nabla_{\feta}\tLoss \left(\tbt,\tft \right) = -\frac{1}{n}\sum_{i=1}^{n} \loss_{2}' (Z_{i}^{T}\tbt, \tft (\x_{i}), Y_{i}) \kont ( \cdot ,\x_{i}) - 2 \lambda_{n} \tft. 
$$
Define the subspace $\mathcal{H}_{1} = \operatorname{Span}\{\kont ( \cdot ,\x_{i}),i=1,2,...,n\}$ within the considered NTK RKHS. It is straightforward to verify that $\tft \in \mathcal{H}_{1}, \forall t \geq 0$, meaning that the evolution of $\tft$ is restricted to this finite-dimensional subspace. 
To analyze the optimization convergence, the following assumption on convexity and smoothness is standard. 
\begin{assumption}[Conditions for optimization]\label{ass:convexloss} 
The loss function $\loss = \loss(\cdot, \cdot, \cdot)$ is convex and non-negative, the gradient $\nabla \loss$ is $\jdc{c1}$-Lipschitz continuous with a constant $\xdc\label{c1}$. 
\end{assumption}

This assumption on the loss function ensures that the gradient flow converges. The convexity guarantees convergence, while the Lipschitz continuity of the gradient ensures stability during optimization, preventing abrupt changes that could hinder convergence. 
Denote the initial loss $L_0 = \tLoss(\beo, \fo)$ and assume that the RKHS global minimizer of \eqref{equ:rkhs} satisfies $ \| \tbw \|_{2}^2 +\|  \tfw \|_{\rkhs}^2 = \tilde{B}^2$. 
The following conclusion characterizes the convergence of the ideal gradient flow for general loss functions. 

\begin{proposition}\label{prop:convergence}
Given the sample \( \{(Y_i, Z_i, \x_i), i = 1, 2, \dots, n\} \) and we consider the optimization problem \eqref{equ:rkhs} in the Hilbert space $\mathcal{H}_{1}$ with the gradient flow method described above, the following results on the optimization convergence rate hold.

(1). When Assumption \ref{ass:convexloss} holds, we have 
$$
   \| \nabla_{\beta} \tLoss (\tbt,\tft ) \|_{2}^2 + \| \nabla_{\feta}\tLoss (\tbt,\tft )\|_{\rkhs}^2 \lesssim \tLoss(\tbt, \tft ) - \tLoss(\tbw, \tfw) \lesssim \frac{\tilde{B}^2\Lo }{\Lo t + \tilde{B}^2}.  
$$

(2). Additionally, if \eqref{equ:rkhs} is $\mu$-strongly convex for a positive number $\mu$, we have 
$$
 \| \nabla_{\beta} \tLoss (\tbt,\tft ) \|_{2}^2 + \| \nabla_{\feta}\tLoss (\tbt,\tft )\|_{\rkhs}^2 \lesssim   \tLoss(\tbt, \tft ) - \tLoss(\tbw, \tfw) \lesssim  e^{- 2\mu t}\Lo.  
$$
\end{proposition}

\begin{remark}\label{rm:trainingtime}
    The above proposition demonstrates that the convergence rate of optimization is sub-linear for convex objects and linear for strongly convex objects; which is standard in optimization theory. 
    Therefore, to obtain the optimizer $(\tbt, \tft )$ that satisfies 
    $$ \tLoss(\tbt, \tft ) - \tLoss(\tbw, \tfw) 
    \quad \text{ and } \quad 
    \|\nabla_{\beta} \tLoss (\tbt,\tft ) \|_{2} + \| \nabla_{\feta}\tLoss (\tbt,\tft)\|_{\rkhs} \lesssim n^{-1}, 
    $$
    we need to train $\stopt$ that 
    $ \stopt \gtrsim n^{2} \tilde{B}^{2}$ when the objective is just convex and $  \stopt \gtrsim \mu^{-1}\log n $ when it is $\mu$ strongly convex. 
    Furthermore, the penalization procedure guarantees $\|  \tfw \|_{\rkhs}^2 \lesssim \Lo/\lambda_n$ and the consistency (or at least boundedness) of $\tbw$ is usually easy to establish in statistical problems. 
    Therefore, without involving neural networks, the upper bounds of $\tilde{B}^2$ for specific problems are typically available. 
\end{remark}

The following theorem demonstrates that, for any given training time, the discrepancy between the neural network and RKHS optimization results can be sufficiently small when the neural network is wide enough.

\begin{theorem}\label{thm:gap} 
Given any positive real number $\xi \in(0,1)$ and the training time $t_{}$. Let the network $m(t,\xi) \geq \poly( n,\lambda^{-1},\Lo,\log(1/\xi), \exp (t_{}) )$ large enough.
Then, with probability at least $1-\xi$ over neural network random initialization, the differences between neural network estimation and RKHS estimation can be bounded by 
$$
\max_{1\leq s \leq t_{} } \left\{ \|\hbs-\tbs\|, \| \hfs - \tfs \|_{L^{\infty}} \right\}  \leq o(n^{-1/2}). 
$$
\end{theorem}

This theorem demonstrates that the gradient flow of DNNs and the NTK RKHS can exhibit significant similarity when the neural network is sufficiently wide. Although the optimization landscape of DNNs is generally complex, in the NTK overparameterized regime, it closely resembles that of the RKHS. This similarity also helps avoid the degeneration of the tangent space as seen in the counterexample in Section \ref{sec:chacon}, with high probability. Additionally, the requirement on the width of DNNs includes an exponent term related to the training time, which is used to bound the accumulated dynamic differences between the two training processes. This factor accounts for the cost of accommodating general, unstructured loss functions, but can be omitted when considering the least squares loss as a special case.

\section{Statistical theory of algorithmic neural M-estimation} \label{sec:str}
In this section, we discuss the statistical properties of the neural semiparametric $M$-estimator, focusing on the convergence rate of the nonparametric part and the asymptotic normality of the parametric component.
Notably, we do not impose common assumptions such as the Lipschitz continuity \citep[e.g.][]{ma2021statistical} or strongly convexity of the loss function, as these conditions may limit the applicability of our results.
Furthermore, because our analysis pertains to the algorithmic optimizer, we do not assume that the nonparametric component is bounded; this introduces additional challenges for the theoretical analysis.
Before discussing details, we introduce some basic assumptions.
Firstly, some basic conditions of the model regularity are summarized below.
\begin{assumption}
(1) The covariate $(Z, \x)$ takes value in a bounded domain with a joint probability density function bounded away from $0$ and $\infty$. 

(2) Conditional on $(Z, \x)$, the second order moment $\mathbb{E}\left[ 
Y^2|Z=\zz, \x=\xx \right]$ is bounded. 

(3) Derivatives and expectations are exchangeable in the sense that
$$
\frac{\partial}{\partial u }
\mathbb{E}\left[ \loss(u_1, u_2,Y)|Z=\zz, \x=\xx \right] = \mathbb{E}\left[ \frac{\partial}{\partial u } \loss(u_1, u_2,Y)|Z=\zz, \x=\xx \right],
$$
for $u=u_1$ and $u_2$. 
\end{assumption}

When establishing the nonparametric convergence, the margin condition \citep{tsybakov2004optimal} and the Bernstein condition \citep{bartlett2006empirical} are often employed to achieve fast rates in statistical and learning analysis. 
Consider a simple nonparametric model \( P_{\feta} \) for example, which depends on a nonparametric parameter \( \feta \) and a corresponding loss function \( \loss_\feta \).
These conditions quantify the identifiability and the curvature of the objective function \( f \mapsto \Pe \loss_f \) at some minimum \( f^* \). 
In the margin condition, \( f^* = \fo \) is the minimizer of the risk over all measurable functions, whereas in the Bernstein condition, \( f^* \) typically minimizes the risk over the candidate function class \( \mathcal{F} \), see \citet{lecue2011interplay} for more discussions. 
As one of specific forms, these conditions may establish relationships between the excess risk \( \Pe\left(\loss_f - \loss_{f^*}\right) \) and the \( L^2 \) norm \( \|f - f^*\|_{L^2}^2 \asymp \Pe\left(\feta(\x) - f^*(\x)\right)^2 \) through inequalities of the form \( \Pe\left(\loss_{\feta} - \loss_{f^*} \right) \gtrsim \|f - f^*\|_{L^2}^{2\kappa} \) with typically $\kappa=1$.
This implies a better concentration and smaller localized sets, 
hence helps for the fast convergence. 
However, to verify the Margin/Bernstein conditions usually requires that $\feta$ is near $\fo$, for example, $\|\feta - \fo\|_{L^\infty}$ is bounded.  
This has limited the application of related results when the boundedness of nonparametric functions does not hold naturally.
Especially in this work, assuming boundedness for the optimizer under gradient flow is particularly unsuitable. 
In semiparametric problems, unlike the finite-dimensional parameters for which consistency is typically easy to establish, $L^{\infty}$ boundedness of the nonparametric estimation is often non-trivial, especially when the true function does not lie within the considered RKHS. 
The following Huberized margin condition holds more easily for unbounded function class.

\begin{assumption}[Huberized margin condition for semiparametric estimation]  \label{ass:subc} 
There is a constant $\xdc\label{c3}>0$ such that for every $\beta \in 
\mathbb{R}^{p}$, $\feta \in \mathcal{F},$
     \begin{equation}\label{equ:submc}
     \Pe \left(\lossbf  - \lossbeofo   \right) \geq \jdc{c3} \frac{\left\|\beta-\beo\right\|^2 + \left\|f-\fo\right\|_{L^2(\x)}^{2}}{ 1 + \left\|\beta-\beo \right\| + \|\feta - \fo \|_{L^\infty} } .
     \end{equation}
\end{assumption}

Typically, when the loss is bounded or globally Lipschitz continuous \citep[e.g.][]{ma2021statistical}, the class complexity is easy to bound, but such assumptions may not be general enough. 
In the proposed Huberized condition, ignoring finite-dimensional parameter $\beta$ for convenience, when $\| \feta-\fo \|_{L^\infty} \lesssim 1$, the right hand is nearly $\left\|\feta-\fo\right\|_{L^2(\x)}^{2}$; while when $\| \feta-\fo \|_{L^\infty} \gtrsim 1$, the right hand is even smaller than $\left\|\feta-\fo\right\|_{L^1(\x)}$. 
If only the local curvature is of concern, then this condition is equivalent to the commonly used margin condition  \( \Pe\left(\loss_{\feta} - \loss_{\fo} \right) \gtrsim \| \feta - \fo \|_{L^2}^{2} \).
As pointed out earlier, it is not proper to assume that the algorithmic neural network optimizer is bounded, and the $L^{\infty}$ consistency or boundedness of nonparametric estimation is usually non-trivial when the true function does not lie within the considered RKHS. 
Moreover, in \eqref{equ:submc} for semiparametric models, $\left\|\beta-\beo \right\|$ in the denominator is for symmetry, which is usually not necessary because the consistency of the estimation of $\beta$ is verifiable in common cases. 

When the risk is strongly convex, the condition is obviously satisfied. 
Intuitively, the proposed assumption only requires that the Hessian matrix is non-singular near the minimum and that the gradient is lower bounded at the far end. 
Strong convexity of the loss function across the entire domain is unnecessary; local strong convexity near the true minimizer is often sufficient. This is not much stricter than Assumption \ref{ass:convexloss}, more examples for illustration will be given in Section \ref{sec:example}.

The next assumption assumes that the true parameters lie within the ideal space. 
\begin{assumption} \label{ass:np1} 
The true parameter $\fo$ belongs to the RKHS $\rkhs$. 
\end{assumption}

This condition is common in traditional statistical works and yields some boundedness to simplify the analysis. However, it does not always hold in a broad context of learning problems and is in fact not necessary to achieve the desired property. Thus, we relax this condition to allow that the true function does not lie within the reproducing space of NTK.

\begin{manualtheorem}{\ref*{ass:np1}$'$}\label{ass:np2} 
The true parameter \(\fo\) does not reside within the RKHS \(\rkhs\), but instead belongs to a Sobolev space \(\sob^{s,2}\) with \(s > d/2\).
\end{manualtheorem}

Here, we assume that the smoothness of \(\fo\) satisfies \(s > d/2\), a condition crucial for establishing the statistical properties under consideration, such as the convergence rate and the attainment of parametric asymptotic normality. This assumption is particularly significant in the context of semiparametric statistics, where achieving a \(\sqrt{n}\)-consistent estimator for finite-dimensional parameters generally necessitates that the non-parametric convergence rate faster than \(n^{-1/4}\) \citep{van2000asymptotic,kosorok2008introduction}.

Now, we introduce some basic concepts and assumptions for the semiparametric model \citep{van2000asymptotic,kosorok2008introduction}.
Given a fixed $\feta$, 
let $\mathcal{G}_{\feta}$ denote the collection of all smooth functional curves $g_b$ that run through $\feta$ at $b=0$, i.e. 
$\{g_b: g_b \in L^2\left( \dom \right)$ for $ b \in(-1,1)$ and 
$\lim _{b \rightarrow 0}\left\|b^{-1}\left(g_b-g_0\right)-h\right\|_{L^2} \rightarrow 0$ for some function $h \}$. 
Then let 
$$\mathcal{T}_{\feta}\mathcal{G}_{\feta} =\left\{h: \dom \to \mathbb{R}, \text{ there is a }  g_b \in \mathcal{G}_{\feta } \text{ such that } \lim _{b \rightarrow 0}\left\|b^{-1}\left(g_b-g_0\right)-h\right\|_{L^2}=0\right\}$$
be all the potential tangent directions for the nuisance parameter. 
Now we set $\overline{\mathcal{T}}_{\feta }\mathcal{G}_{\feta}$ be the closed linear span of $\mathcal{T}_{\feta }\mathcal{G}_{\feta}$ under linear combinations.

For any $h$ in $\mathcal{T}_{\feta}$, denote
$$
S_1(\beta, \feta) = 
\frac{\partial}{\partial \beta} \lossbf  \text { and } S_2(\beta, \feta)[h]=\left.\frac{\partial}{\partial b}\right|_{b=0} \loss_{\beta, g_b}, 
$$
where $g_b$ satisfies $\left.(\partial /(\partial b)) g_b\right|_{b=0}=h$. 
For convenience, the tangent set for $\feta$ at $(\beta,\feta)$ is defined as $\Lambda_{\beta,\feta}:=\{ S_2(\beta, \feta)[h]: h\in \mathcal{T}_{\feta } \}$.  
Further, we also define 
$$
\begin{gathered}
S_{11}(\beta, \feta)=\frac{\partial}{\partial \beta} S_1(\beta, \feta), 
\quad S_{12}(\beta, \feta)[h]=\left.\frac{\partial}{\partial b}\right|_{b=0} S_1(\beta, g_b), \\
S_{21}(\beta, \feta)[h]=\frac{\partial}{\partial \beta} S_2(\beta, \feta)[h]\text { and } 
S_{22}(\beta, \feta)\left[h_1\right]\left[h\right]=\left.\frac{\partial}{\partial b}\right|_{b=0} S_2\left(\beta, g_b\right)\left[h_1\right],
\end{gathered}
$$
where $h_1$ is another function in $\mathcal{T}_{\feta}$.

As is well known, in the special case when the loss is a negative log-likelihood, the efficient score is essentially the projection of the score function $S_1$ onto the orthonormal complement of the tangent set $\Lambda_{\beta,\feta}$. Denote $\Pi_{\theta, \feta}$ as the orthogonal projection onto the closure of the linear span of $\Lambda_{\beta,\feta}$. The efficient score function for $\beta$ at \((\beo, \fo)\) is then given by
$$
\seff(\beo, \fo)=S_1(\beo, \fo)-\Pi_{\beo, \fo} S_1(\beo, \fo).
$$ 
Denote $$S_2(\beta, \feta)[\ddh]=\left(S_2(\beta, \feta)\left[h_1\right],S_2(\beta, \feta)\left[h_2\right], \ldots, S_2(\beta, \feta)\left[h_p\right]\right),$$ where $\ddh=\left(h_1,h_2, \ldots, h_p\right) \in \mathcal{T}_{\feta}^p$. 
Then define $S_{12}\left[\ddh \right]$, $S_{21}\left[\ddh\right]$  and $S_{22}\left[\ddh_1\right]\left[\ddh_2\right]$ accordingly. If there is $\tddh=(\thh_1, \thh_2, \ldots, \thh_p^*) \in \mathcal{T}_{\feta}^p\mathcal{G}_{\feta}$, such that for any $\ddh=\left(h_1, \ldots, h_k\right) \in \mathcal{T}_{\feta}^p\mathcal{G}_{\feta}$,
\begin{equation}\label{equ:lfd}
    \Pe\left(S_{12}\left(\beo, \fo\right)[\ddh]-S_{22}\left(\beo, \fo\right)[\tddh ][\ddh]\right)=0. 
\end{equation}
Then we can also determine the efficient score for likelihood estimation as 
$$
\seff(\beo, \fo) = S_{1}\left(\beo, \fo\right)-S_{2}\left(\beo, \fo\right)[\tddh ]. 
$$
Similar to the negative log-likelihood, the following assumptions are made for the general loss function.

\begin{assumption}\label{ass:semi} 
(1). (Positive information) There exists $\tddh=(\thh_1, \ldots, \thh_p) \in \mathcal{T}_{\feta}^p\mathcal{G}_{\feta}$ such that \eqref{equ:lfd} holds for any $\ddh=\left(h_1, \ldots, h_k\right) \in \mathcal{T}_{\feta}^p\mathcal{G}_{\feta}$. Further assume that 
$$
\thh_i^* \in \sob^{s,2}(\dom), i=1,2,...,p, 
$$
where $s>d/2$. Given $\tddh$, the matrix 
$$
\begin{aligned}
A = \Pe\left(S_{11}\left(\beo, \fo\right)-S_{21}\left(\beo, \fo\right) [\tddh ]\right)
\end{aligned}
$$
is nonsingular. 

(2). (Smoothness of the model) For all possible parameters $(\beta, \feta)$ satisfying $\{(\beta, \feta):\left\|\beta-\beo\right\|$$+ \left\|\feta-\fo \right\|_{L^2(\x)} \lesssim \delta_n \}$ with $ \delta_n = o(n^{-1/4})$, 
$$
\begin{aligned}
    & \Pe \Big\{(\seff(\beta, \feta)-\seff\left(\beo, \fo\right))-\left(S_{11}\left(\beo, \fo\right)-S_{21}\left(\beo, \fo\right)[\tddh ]\right)\left(\beta-\beo\right) \\
&\quad \quad \quad \quad \quad \quad \quad \quad \quad \quad \quad \left.-\left(S_{12}\left(\beo, \fo\right)\left[\feta-\fo \right]-S_{22}\left(\beo, \fo\right)[\tddh ]\left[\feta-\fo\right]\right)\right\}\\
&\quad \quad \quad =o\left(\left\|\beta-\beo\right\|\right)+O(\left\|\feta-\fo \right\|_{L^2}^{2}) .
\end{aligned}
$$

\end{assumption}

In the above assumptions, the first one guarantees the existence and smoothness of directions $\tddh$, which corresponds to the least favorable direction when the loss is the negative log-likelihood. 
This direction can be determined by solving the equations in \eqref{equ:lfd}.
The second condition can usually be obtained by Taylor expansion. 
Under the above assumptions, we have the following general conclusion.

\begin{theorem}\label{thm:main}
Consider the proposed semiparametric neural $M$-estimation, assume that the estimation $(\hbt,\hft)$ is trained from \eqref{equ:011},\eqref{equ:012},\eqref{equ:013} with training time $\stopt$ in Remark \ref{rm:trainingtime} and network width $m(\stopt,\xi)$ in Theorem \ref{thm:gap} with $\xi \in(0,1)$. 
Then, with probability at least $1-\xi$ over neural network random initialization, the following results holds: 
If $\hbts$ and $\Pe \loss_{\hbts,\hfts}$ are bounded, and Assumption \ref{ass:convexloss}, \ref{ass:subc}, \ref{ass:np1},  \ref{ass:semi} hold, 
setting $\lambda \asymp n^{-(d+1)/(2d+1)}$, 
we have 
\begin{equation}\label{equ:nprate}
    \| \hfts - \fo \|_{L^2}^2 = O_p\left( n^{-(d+1)/(2d+1)} \log n \right).  
\end{equation}
If Assumption \ref{ass:np1} is replaced by Assumption \ref{ass:np2}, setting $\lambda \asymp n^{-(d+1)/(2s+d)}$, the nonparametric rate becomes 
\begin{equation}\label{equ:npratep}
    \| \hfts - \fo \|_{L^2}^2 = O_p\left( n^{-2s/(2s+d)} \log n \right). 
\end{equation}
In both cases, we have the asymptotic normality
\begin{equation}\label{equ:normal}
\begin{aligned}
&\sqrt{n}\left(\hbts-\beo\right)=
n^{1/2}  A^{-1}  \Pn \seff(\beo, \fo) +o_p(1) 
\xrightarrow{d} N(0, \Sigma). 
\end{aligned}
\end{equation}
where $\Sigma=A^{-1} (\Pe \{ \seff(\beo, \fo)\seff(\beo, \fo)^{T} \}) (A^{-1})^T $. 
\end{theorem}

The boundedness condition of $\hbts$ and $\Pe \loss_{\hbts,\hfts}$ is essentially a boundedness requirement for the parameter part, not for the nonparametric $\feta$. 
For nonparametric problems without the finite-dimensional parameter $\beta$, this condition can be removed under the same proof. 
For semiparametric models, this condition is often verifiable in specific examples, as in the applications in the next section. 
Let $\so = (d+1)/2$ denote the smoothness of the RKHS corresponding to the NTK. The convergence rate in \eqref{equ:nprate} can then be expressed as nearly $n^{-2\so/(2\so + d)}$ with the tuning parameter $\lambda \asymp n^{-2\so/(2s+d)}$. This smoothness is essentially determined by the ReLU activation function, and the rate is minimax optimal for many statistics models with Assumption \ref{ass:np1}. When the true function does not lie within the considered ReLU NTK RKHS, as stated in Assumption \ref{ass:np2}, the rate in \eqref{equ:npratep} remains minimax optimal. Furthermore, Assumption \ref{ass:np2} ensures that the nonparametric rate is faster than $n^{-1/4}$.
Therefore, regardless of whether the true function resides in the considered RKHS, the estimator is shown to be asymptotically normal.

An estimator is semiparametric efficient if its information (i.e., the inverse of its variance) equals the supremum of the information across all parametric submodels. The submodel that attains this supremum is typically called the least favorable or hardest submodel. By the above result, when the loss function is the negative log-likelihood and the model class contains the hardest submodel, the proposed semiparametric estimator is efficient. Conversely, if the loss function differs, suggesting model misspecification, the parametric estimator remains $\sqrt{n}$-consistent and asymptotically normal.

We conclude this section with the following remarks. To better align with practice, we consider overparameterized neural networks and general loss functions, allowing the true function to lie outside the desired RKHS. Distinct from existing literature, we refrain from assuming that the output of the optimized network is bounded and study the algorithmic solution, which presents significant challenges for statistical analysis and inspires a new yet weaker margin condition. Then, by combining the peeling argument with the interpolation bound, we precisely characterize the entropy and establish the nonparametric convergence rate. Based on this, we obtain the asymptotic normality of the parametric estimators, with high probability over the random initialization, avoiding the degeneration of the tangent space and addressing the key question posed in Section \ref{sec:chacon}. Lastly, we acknowledge the potential for broader applications of the proposed framework. This essentially provides a fundamental approach to addressing statistical problems that involve characterizing the tangent spaces of the nonparametric component, which has been a challenge that frequently arises in the intersection of deep learning and statistical inference.

\section{Examples}\label{sec:example}
In this section, we introduce two common examples to illustrate the proposed semiparametric neural $M$-estimator and the established theoretical framework.

\subsection{Regression}
Now we consider the regression model 
\begin{align}
\label{eq:partial_linear_regression}
   Y_i = \fo(X_{i}) + Z_{i}^{T}\beo + \epsilon_i,  
\end{align}
where $Y_i$ represents the response variable, $X_i$ and $Z_i$ are bounded covariates with densities bounded away from $0$ and $\infty$, and $\epsilon_i$ denotes the independent measurement noise with mean zero, finite variance and a symmetric distribution. For estimation, we define the loss function as  
$\loss=\loss(Y_i- \feta(X_{i}) - Z_{i}^{T}\beta)$, where $l$ is a univariate function. 

The following condition is for the covariates distribution and loss function. 

\begin{assumption}\label{ass:lossregression}
~\par
(1) The matrix $\mathbb{E}\left[(Z-\mathbb{E}[Z|\x]) (Z-\mathbb{E}[Z|\x])^T \right]$ is nonsingular. 

(2) The univariate function $\loss$ is a nonnegative, even and convex function with Lipschitz continuous derivative $\loss'$ and $\loss(0) = \loss'(0) = 0$. 
Additionally, the pointwise risk $R(s) = \mathbb{E}_{\epsilon}[\loss(s+\epsilon)]$ is convex and 
$R''(s)\geq \jdc{c5}>0$ for 
$s: |s|\leq b$ with some $b$ large enough 
and 
a positive number $\xdc\label{c5}$. 
\end{assumption}

In this example, symmetry of the noise and loss functions is assumed to ensure that \( (\beo, \fo) \) is the minimizer of the risk \( \mathbb{E}_P \loss_{\beta, \feta} \), thereby simplifying the analysis. 
Condition (1) in Assumption \ref{ass:lossregression} is standard \citep{kosorok2008introduction} to ensure
$$\mathbb{E}\Big[ \big( Z^{T}\beta + \feta(\x) - Z^{T}\beo - \fo(\x)  \big)^2 \Big] \gtrsim \|\beta - \beo\|^2 + \| \feta - \fo \|_{L^2}^2, $$ 
in partially linear models. 
Furthermore, the second condition assumes that the risk function is strongly convex within a local neighborhood. Many commonly used loss functions, such as least squares loss and Huber loss, satisfy this assumption. Assumption \ref{ass:lossregression} guarantee that Assumption \ref{ass:subc} holds, as shown in Lemma \ref{lem:explm} in the Supplementary Material \citep{supplyment}.

\begin{theorem}\label{cor:reg}
Consider the proposed semiparametric neural $M$-estimation for model \eqref{eq:partial_linear_regression} and set the hyperparameters as in Theorem \ref{thm:main}. Under Assumption \ref{ass:np1}, \ref{ass:semi}(1), \ref{ass:lossregression}, 
setting $\lambda \asymp n^{-(d+1)/(2d+1)}$, with probability at least $1-\xi$ over neural network random initialization, 
we have 
$$
\| \hfts - \fo \|_{L^2}^2 = O_p\left( n^{-(d+1)/(2d+1)} \log n \right). 
$$
If Assumption \ref{ass:np1} is replaced by Assumption \ref{ass:np2}, setting $\lambda \asymp n^{-(d+1)/(2s+d)}$, the nonparametric rate becomes 
$$
    \| \hfts - \fo \|_{L^2}^2 = O_p\left( n^{-2s/(2s+d)} \log n \right).
$$ 
In both cases, we have the asymptotic normality
\begin{equation}\label{equ:normal3}
\begin{aligned}
&\sqrt{n}\left(\hbts-\beo\right)=
n^{1/2}  A^{-1}  \Pn \seff(\beo, \fo) +o_p(1) 
\xrightarrow{d} N(0, \Sigma), 
\end{aligned}
\end{equation}
where $\Sigma=A^{-1} (\Pe \{ \seff(\beo, \fo)\seff(\beo, \fo)^{T} \}) (A^{-1})^T $. 
\end{theorem}

Some previous works studied related challenges in the least squares nonparametric regression problem \citep{mendelson2010regularization, steinwart2009optimal}. 
In the above theorem, whether the true function belongs to the considered RKHS, we establish the nonparametric optimal convergence rate and parameter asymptotic normality, without assuming that candidate functions are bounded and allowing unbounded responses. 


\subsection{Classification}
In this subsection, we consider the following binary classification problem.
Suppose that we can observe independently and identically distributed random sample $\{(Y_i,Z_i,\x_i),i=1,2,\cdots,n\}$, where $Y_i \in \{0,1\}$ denotes binary response, and $\x_{i}'s$ and $Z_{i}'s$ are bounded covariates with densities bounded away from $0$ and $\infty$. We assume that $Y$ follows a Bernoulli distribution determined by:
\begin{equation}
    \label{eq:partial_linear_classification}
    \ppp(Y = 1) = \phi\left(\fo(X) + Z^{T}\beo\right)  \quad\text{   and  }\quad  \ppp(Y = 0) = 1-\phi\left(\fo(X) + Z^{T}\beo\right),
\end{equation}
where $\phi: \mathbb{R} \mapsto [0,1]$ is a continuously differentiable monotone link function. 
Hence the loss function can be taken as the negative log-likelihood as
\begin{equation}\label{equ:lossclassif}
    \loss(Z^{T}\beta, \feta(X), Y) = -Y \log \phi\left(\feta(X) + Z^{T}\beta\right) - (1-Y )\log \left(1-\phi\left(\feta(X) + Z^{T}\beta\right)\right). 
\end{equation}
Common choices for $\phi$ include the sigmoid function $\phi(t)=(1+e^{-t})^{-1}$ and the cumulative normal distribution function, corresponding to the logistic and probit models, respectively. 
In theoretical analysis, we allow for model misspecification, assuming that the data-generating process follows: 
\begin{equation}\label{eq:partial_linear_classification_mis}
    \ppp(Y = 1) = \psi\left(\feta_{1}(X) + Z^{T}\beta_{1}\right)  \quad\text{   and  }\quad  \ppp(Y = 0) = 1-\psi\left(\feta_{1}(X) + Z^{T}\beta_{1}\right),
\end{equation}
where $\psi \neq \phi$ is a different link function, and some $\beta_{1} \in \mathbb{R}^{p},\feta_{1} \in L^2(\dom)$. 
Despite this potential misspecification, we continue to use the working link function $\phi$ and loss function as specified in \eqref{equ:lossclassif}. 
From a variational perspective, the minimizer $(\beo,\fo)$ of the $\Pe\loss_{\beta,\feta}$ can also be interpreted as minimizer of the Kullback-Leibler divergence 
$$
\Pe\left[\frac{\log P_{\beta, \feta}}{\log P_0}\right]=
\Pe\left[\frac{Y \log \phi\left(\feta(X ) + Z ^{T}\beta\right) + (1-Y )\log \left(1-\phi\left(\feta(X ) + Z^{T}\beta\right)\right)}
{Y \log \psi\left(\feta_{1}(X ) + Z ^{T}\beta_{1}\right) + (1-Y )\log \left(1-\psi\left(\feta_{1}(X ) + Z^{T}\beta_{1}\right)\right)}\right].
$$
Thus, the parameters $(\beo, \fo)$ retain interpretability and are still parameters in terms of the underlying distribution. 
Solving the equation \eqref{equ:lfd}, some calculation implies
$$ 
\tddh= \frac{
\mathbb{E}\left[ Z\left(Y\left(\log \phi \right)^{''} + (1-Y)\left(
\left(\log ( 1-\phi) \right)^{''}\right)\right) |\x \right]
}{
\mathbb{E}\left[Y\left(\log \phi \right)^{''} + (1-Y)\left(
\left(\log ( 1-\phi) \right)^{''}\right) |\x \right]
}. 
$$

The following conditions are standard for 
the link function. 

\begin{assumption}\label{ass:clas}
%
The link function $\phi$ is Lipschitz continuous, monotone and continuously differentiable; $-\log\phi(t)$ and $-\log(1-\phi(t))$ are both convex with positive, continuous and bounded second-order derivative.
If the model is misspecified, we make the same assumptions for $\psi$ and assume that the underlying $\beta_1, \feta_1$ are bounded in infinity norm. 

\end{assumption}

This condition ensures that the loss function and distribution satisfy Assumptions \ref{ass:convexloss} and \ref{ass:subc}. Commonly used link functions, such as the logistic and probit functions, naturally satisfy this condition, thereby making these assumptions applicable in practice.

\begin{theorem}\label{cor:classi}
Consider the proposed semiparametric neural $M$-estimation for model \eqref{eq:partial_linear_classification} or \eqref{eq:partial_linear_classification_mis} and set the hyperparameters as in Theorem \ref{thm:main}. Under Assumption \ref{ass:np1}, \ref{ass:semi}(1), \ref{ass:lossregression}(1) and \ref{ass:clas}, 
setting $\lambda \asymp n^{-(d+1)/(2d+1)}$, with probability at least $1-\xi$ over neural network random initialization, 
we have 
$$
\| \hfts - \fo \|_{L^2}^2 = O_p\left( n^{-(d+1)/(2d+1)} \log^2 n \right).
$$
If Assumption \ref{ass:np1} is replaced by Assumption \ref{ass:np2}, setting $\lambda \asymp n^{-(d+1)/(2s+d)}$, the nonparametric rate becomes 
$$
    \| \hfts - \fo \|_{L^2}^2 = O_p\left( n^{-2s/(2s+d)} \log^2 n \right).
$$ 
In both cases, we have the asymptotic normality
\begin{equation}\label{equ:normal2}
\begin{aligned}
&\sqrt{n}\left(\hbts-\beo\right)=
n^{1/2}  A^{-1}  \Pn \seff(\beo, \fo) +o_p(1) 
\xrightarrow{d} N(0, \Sigma), 
\end{aligned}
\end{equation}
where $\Sigma=A^{-1} (\Pe \{ \seff(\beo, \fo)\seff(\beo, \fo)^{T} \}) (A^{-1})^T $. 
\end{theorem}

For the partially linear classification problem, the proposed neural $M$-estimation achieves the minimax nonparametric rate as well as $\sqrt{n}$-consistency, with high probability. This performance is attributed to the representational capacity of overparameterized deep neural networks and their tangent space. When the model is correctly specified, i.e. the loss function corresponds to the negative log-likelihood, the resulting parametric estimator attains efficiency.

\section{Numerical studies} \label{sec:num}

This section demonstrates the practical advantages of the proposed semiparametric neural $M$-estimator through simulation studies on both regression and classification models described in Section \ref{sec:example}. We use an overparameterized fully connected ReLU neural network with depth $L=5$ and width $m=1000$. For training, the stochastic gradient descent optimizer in PyTorch is employed,  with a learning rate of $0.001$ and a total of $1000$ epochs.
For a comprehensive evaluation, we compare our method with four baseline approaches that also estimate \(\beo\) and \(\fo(x)\) via (penalized) $M$-estimation.
The first baseline is a regression spline estimator using B-splines, with uniformly spaced knots over the domain \([0,1]^d\). The second is the RKHS method, where the RKHS norm of the nonparametric component serves as the penalty, and the Laplacian kernel \(K_h(x_1, x_2) = \exp\{-\|x_1 - x_2\|/h\}\) is used. The third is a local linear estimator with the Epanechnikov kernel \(k_h(x_1,x_2) = (1 - \|x_1 - x_2\|^2/h)/(2^d(1 - 1/(3h)))\).
The last one also employs the fully connected ReLU neural network with depth $L=5$, and the width serves as a hyperparameter, referred to as ``underparameterized'' to be distinguished from the overparameterized regime.   

To select the optimal hyperparameters for all methods, we split the full data into a training set (80\%) and a validation set (20\%) and choose the hyperparameters that minimize validation loss, including the tuning parameters for our method and the second baseline, the number of basis functions for the first baseline, the bandwidth for the third baseline, and the width of the underparameterized neural network. 
Specifically, the hyperparameters of the regression spline method and the underparameterized neural network method are chosen from a set that ensures the total number of parameters is less than the sample size.


We generate \(\{Z_i=(Z_{1i},Z_{2i})^T\}\), \(i=1, 2, \dots, n\) from a uniform distribution over the interval $[0,1]^2$. Then, \(\{X_i=(X_{1i},\cdots,X_{di})^T\}\) is generated using the formula $X_{ji}=0.9W_{ji}+0.05(Z_{1i}+Z_{2i}),1\le j\le d$, where $W_{ji}$ is sampled from a uniform distribution over the interval $[0,1]$. 
The true finite-dimensional parameter vector is set as \(\beta_0 = (1, 0.75)^T\).
For the nonparametric part \(\fo(x)\), we consider four cases with different dimensions and function forms. Here, Case 1 and Case 3 correspond to five-dimensional examples, while Case 2 and Case 4 represent their respective extensions to ten-dimension; which have been studied in the simulation of \citet{deepcox} and \citet{yan2023nonparametricregressionrepeatedmeasurements}.
\begin{align*}
    &\operatorname{Case\ 1:} d = 5,\ f_0(x) = 5\left(x_1^2x_2^3+\log(1+x_3)+\sqrt{1+x_4x_5} +\exp(x_5/2)\right),\\
    &\operatorname{Case\ 2:} d = 10,\ f_0(x) =\frac{5}{2}\left(x_1^2x_2^3+\log(1+x_3)+\sqrt{1+x_4x_5} +\exp(x_5/2)\right. \\&\hspace{10em}+\left.x_6^2x_7^3+\log(1+x_8)+\sqrt{1+x_9x_{10}} +\exp(x_{10}/2)\right),\\
      &\operatorname{Case\ 3:} d = 5,\ f_0(x) = 5\sin\left(\frac{6\pi}{d(d+1)}\sum_{l=1}^5 lx_l\right), \\
    &\operatorname{Case\ 4:} d = 10,\ f_0(x) = 5\sin\left(\frac{6\pi}{d(d+1)}\sum_{l=1}^{10} lx_l\right).
\end{align*}

\begin{table}[tb]
\centering
\centerfloat
    \caption{The mean square error $(\times 10^{-1})$ for $\hat{\beta}$ and $\hat{\feta}$ of our method and baselines for the regression model.}
	\label{tab:simulation1}
	\begin{tabular}{cc| ccccc}
 \toprule[1pt]
  \multicolumn{2}{c|}{Setting} & Proposed & Spline & RKHS & Local Linear & Underpara \\
    Case & $n$ &  \multicolumn{5}{c}{MSE for $\hat{\beta}$} \\
    \midrule[1pt]
    \multirow{3}{*}{Case 1}
    & 500 &  0.2613  &  0.4987  &  0.6905  & 0.4148 & 2.0760 \\
    &1000 &  0.1140  &  0.1411  &  0.2511  & 0.1696  & 1.9988 \\
    &2000 &  0.0623  &  0.0513  &  0.0793  & 0.1607 & 1.9965 \\
    \multirow{3}{*}{Case 2}
    & 500 &  0.2429  &  --  &  3.7362  & 132.79  & 3.1735 \\
    &1000 &  0.1132  &  -- &   0.7956  & 20.478  & 3.0583\\
    &2000 &  0.0494  &  --  &  0.1922   & 2.1419  & 3.0356 \\
    \multirow{3}{*}{Case 3}
    & 500 &  0.2316  &  0.6546  &  0.3940  & 116.13  & 7.7303\\
    &1000 &  0.1066  &  0.1801  &  0.1895  & 78.535  & 8.1188 \\
    &2000 &  0.0496  &  0.0755  &  0.0662  &  119.39 & 8.1298\\
    \multirow{3}{*}{Case 4}
    & 500 &  0.2598  &  --  & 0.4656   & 28.555 & 12.192 \\
    &1000 &  0.1098  &  --  & 0.1908   &  88.956 & 12.195 \\
    &2000 &  0.0457  &  --  & 0.0836   & 114.98  & 12.424\\
        \hline
   && \multicolumn{5}{c}{MSE for $\hat{f}$} \\
     \hline
\multirow{3}{*}{Case 1}
    & 500 &  0.8625  &  46.266  & 6.6327   &  2.5206 & 2.8144 \\
    &1000 &  0.7015  &  3.4216  &  3.0720  & 2.4176  & 2.6657 \\
    &2000 &  0.6179  &  0.7906  & 1.5700   & 2.3582  & 2.6089\\
    \multirow{3}{*}{Case 2}
    & 500 &  0.7824  &  --  &  17.558  &  67.927& 2.8644 \\
    &1000 &  0.6084  &  --  &  7.4870  & 11.378  & 2.7192 \\
    &2000 &  0.5093  &  --  &  3.5917  &  2.1836 & 2.6467 \\
    \multirow{3}{*}{Case 3}
    & 500 &  0.6153  &  67.462  &  4.6669  &  89.284 & 40.214 \\
    &1000 &  0.4135  &  7.1053  &  2.3364  &  84.843 & 50.257 \\
    &2000 &  0.2956 &  2.6464   &  1.2638  &  89.256  &  68.170\\
    \multirow{3}{*}{Case 4}
    & 500 &  0.9107  &  --  &  10.753   &  75.208 & 51.472 \\
    &1000 &  0.4489  &  --  &  5.5210  & 84.909  & 51.209 \\
    &2000 &  0.2606  &  --  &  3.0471  & 93.994  & 51.154 \\
         \bottomrule[1pt]
	\end{tabular}
  \end{table}

\begin{table}[tb]
\centering
\centerfloat
    \caption{The mean square error ($\times 10^{-1}$) for $\hat{\beta}$ and $\hat{\feta}$ of our method and baselines for the classification model \eqref{eq:partial_linear_classification}.}
	\label{tab:simulation2}
		\begin{tabular}{cc| ccccc}
 \toprule[1pt]
  \multicolumn{2}{c|}{Setting} & Proposed & Spline & RKHS & Local Linear & Underpara \\
    Case & $n$ &  \multicolumn{5}{c}{MSE for $\hat{\beta}$} \\
    \midrule[1pt]
    \multirow{3}{*}{Case 1}
    & 500 &  0.8821  &  13.561  &  8.8990  & 6.5733  &  2.4930\\
    &1000 &  0.4411  & 12.637   &  2.9707  &  2.3129 & 1.7517 \\
    &2000 &  0.2133 &  12.583  &  1.5815  &  1.8543 & 1.6526 \\
    \multirow{3}{*}{Case 2}
    & 500 &  2.1877  &  --  &  8.1133  & 12.787  &  4.4772\\
    &1000 &  0.5395  &  --  &  2.6608  &  11.939 &  4.1563\\
    &2000 &  0.2177  &  --  &  1.5869  &  11.583 &  4.1018\\
    \multirow{3}{*}{Case 3}
    & 500 &  1.0119  &  21.165  &  22.105  & 13.400  & 4.4189 \\
    &1000 &  0.5900 &  19.280  &  4.2309  & 11.161  & 4.7265 \\
    &2000 &  0.3318 &  17.656  &  1.1527  & 8.0244  & 4.5130 \\
    \multirow{3}{*}{Case 4}
    & 500 &  1.6781  &  --  &  5.5813  &  20.256 & 3.8216 \\
    &1000 &  0.6652  &  --  &  2.3483  & 19.570  & 3.7650 \\
    &2000 &  0.2812  &  --  &  2.5641  & 19.158  & 3.6724\\
        \hline
   && \multicolumn{5}{c}{MSE for $\hat{f}$} \\
     \hline
\multirow{3}{*}{Case 1}
    & 500 &  4.5624 &  39.087  &  9.0062  & 61.677  &  6.9442\\
    &1000 &  2.5099 & 39.043   &  5.0402  & 59.122  & 5.6270 \\
    &2000 &  1.7524  &  38.999  &  3.9864  & 58.263  & 5.2535 \\
    \multirow{3}{*}{Case 2}
    & 500 &  22.417  &  --  &  7.5404  &  29.051 & 11.793 \\
    &1000 &  5.0438  &  --  &  4.0622  &  28.991 & 11.529 \\
    &2000 &  1.8980  &  --  &  3.1317  &  29.023 & 10.751 \\
    \multirow{3}{*}{Case 3}
    & 500 &  7.0489  &  96.863  &  29.990  &  92.969 & 85.572 \\
    &1000 &  3.5228  &  96.762  &  14.225  & 91.686  & 86.227 \\
    &2000 &  2.2837  &   96.850 &  11.259  & 90.390  & 86.411 \\
    \multirow{3}{*}{Case 4}
    & 500 &  17.616  &  --  &  38.825  & 51.957  & 46.156\\
    &1000 &  5.9230  &  --  &  35.755  & 51.703 &  45.684\\
    &2000 &  3.6942  &  --  &  35.791  & 51.598  & 45.623\\
         \bottomrule[1pt]
	\end{tabular}
  \end{table}

In the regression setting, the responses \(Y_i\) are generated by the model \(Y_i = f_0(X_i) + Z_i^T \beta_0 + \varepsilon_i\), where \(\varepsilon_i\) are i.i.d. normal noise with zero mean and standard deviation \(\sigma = 0.5 \). 
For the classification model, the responses \(Y_i\) are drawn from the Bernoulli distribution with probability \(\mathbb{P}(Y_i = 1) = \phi(f_0(X_i) + Z_i^T \beta_0 - \mathbb{E}_X[f_0(X)])\), where \(\phi(x)\) is the logistic function \(\phi(x) = 1/(1 + e^{-x})\). Simulations are performed for three different sample sizes \(n = 500, 1000, 2000\), with $200$ repetitions for each case and method. 
The mean squared error (MSE) of the parametric and nonparametric components is computed, respectively, to evaluate the performance of each method. 
Due to computational limitations, the spline method can only handle Case 1 and Case 3 with $5$-dimensional nonparametric functions.
Tables \ref{tab:simulation1} and \ref{tab:simulation2} report the average MSEs for regression and classification examples, respectively. 
These results demonstrate that our overparameterized neural $M$-estimation approach outperforms the three traditional statistical methods, as well as the underparameterized neural network (i.e., properly tuned with fewer parameters than the sample size). Specifically, due to the high dimensionality of the nonparametric component, the curse of dimensionality becomes significant. Consequently, the underparameterized neural network, spline and the local linear estimator suffer from insufficient learning capacity. 
In all four cases, our method yields the lowest MSE, demonstrating its superior ability.
Lastly, we would like to point out that, in the existing literature \citep[e.g.][]{deepcox,fan2022factor,wang2024deep,yan2023nonparametricregressionrepeatedmeasurements}, regardless of the theoretical requirements imposed on the network size, their numerical experiments have in fact employed  overparameterized neural networks to achieve favorable performance, which also provide certain support for our proposed method and theory. 

Given that the estimation of the parametric component satisfies asymptotic normality, it is possible to perform valid statistical inference. To assess it, we simulate the empirical coverage probabilities of the corresponding estimated confidence intervals.
In the regression model, the variance of $\hat{\beta}$ is estimated as
\begin{align*}
    \hat{\Sigma} = \left[\frac{1}{n}\sum_{i=1}^n (Y_i - \hat{f}(X_i) - Z_i^T\hat{\beta})^2\right]\left[\min_{h \in \mathcal{F}}\frac{1}{n}\sum_{i=1}^n(Z_i - h(X_i))(Z_i - h(X_i))^T\right]^{-1}.
\end{align*}
In the classification model, the variance of $\hat{\beta}$ is estimated as 
\begin{align*}
    \hat{\Sigma} 
    = \left[\min_{h \in \mathcal{F} }\frac{1}{n}\sum_{i=1}^n\phi(\hat{f}(X_i)+Z_i^T\hat{\beta})(1-\phi(\hat{f}(X_i)+Z_i^T\hat{\beta}))(Z_i - h(X_i))(Z_i -h(X_i))^T\right]^{-1},
\end{align*}
where $\mathcal{F}$ is also a neural network function class. 
The coverage rate is defined as the proportion of repeated experiments for which the true parameter falls within the confidence interval. Tables \ref{tab:simulation3} and \ref{tab:simulation4} report the coverage of the $95\%$ confidence intervals constructed by the proposed method based on $500$ repeated experiments, for each case in both regression and classification models.
Generally, the coverage rate is near $0.95$ for the proposed
overparameterized neural $M$-estimation method, especially when the sample size $n$ is large. 
This supports the potential usefulness of our proposed approach for semiparametric inference.

 \begin{table}[tb]
 \centering
\centerfloat
    \caption{The coverage probability for constructed $95\%$ confidence interval \\ for $\beta = (\beta_1,\beta_2)$ in the regression model. }
	\label{tab:simulation3}
	\begin{tabular}{c| ccc|ccc}
 \toprule[1pt]
  Model & \multicolumn{3}{c}{The coverage rate for $\beta_1$} &  \multicolumn{3}{c}{The coverage rate for $\beta_2$}\\
  \hline
 Setting &  $n=500$ &$n=1000$ &$n=2000$&$n=500$ &$n=1000$&$n=2000$ \\
    \midrule[1pt]
    Case 1 & 0.938 & 0.958 & 0.946 & 0.944 & 0.952 & 0.940 \\
    Case 2 &0.966  & 0.958 & 0.950 & 0.960 & 0.950 & 0.952 \\
    Case 3 & 0.946 & 0.970 & 0.948 & 0.968 & 0.952 & 0.950 \\
    Case 4 & 0.986 & 0.974 & 0.938 & 0.988 & 0.968 & 0.948 \\
    \bottomrule[1pt]
	\end{tabular}
  \end{table}

 \begin{table}[tb]
 \centering
\centerfloat
    \caption{The coverage probability for constructed $95\%$ confidence interval\\ for $\beta = (\beta_1,\beta_2)$ in the classification model. }
	\label{tab:simulation4}
	\begin{tabular}{c| ccc|ccc}
 \toprule[1pt]
  Model & \multicolumn{3}{c}{The coverage rate for $\beta_1$} &  \multicolumn{3}{c}{The coverage rate for $\beta_2$}\\
  \hline
 Setting &  $n=500$ &$n=1000$ &$n=2000$&$n=500$ &$n=1000$&$n=2000$ \\
    \midrule[1pt]
    Case 1 & 0.964 & 0.968 & 0.950 & 0.966 & 0.954 & 0.938 \\
    Case 2 & 0.936 & 0.950 & 0.952 & 0.950 & 0.962 & 0.940 \\
    Case 3 & 0.942 & 0.926 & 0.932 & 0.972 & 0.942 & 0.930 \\
    Case 4 & 0.940 & 0.946 & 0.942 & 0.970 & 0.928 & 0.962 \\
    \bottomrule[1pt]
	\end{tabular}
  \end{table}